\documentclass[11pt]{article}

\usepackage{amsmath,amssymb,amsthm}
\usepackage[hidelinks]{hyperref}
\usepackage{geometry}
\numberwithin{equation}{section}

\newtheorem{theorem}{Theorem}[section]
\newtheorem{lemma}[theorem]{Lemma}
\newtheorem{proposition}[theorem]{Proposition}
\newtheorem{corollary}[theorem]{Corollary}
\newtheorem{remark}[theorem]{Remark}

\newcommand{\T}{\mathbb T}
\newcommand{\Td}{\mathbb T^d}
\newcommand{\N}{\mathbb N}

\newcommand{\supp}{\operatorname{supp}}

\title{Sampling recovery for intersections of Lipschitz
classes in the low-smoothness regime}

\author{
Yurii Kolomoitsev\\
\small Institute for Numerical and Applied Mathematics\\
\small University of G\"ottingen, Germany\\
\small \texttt{kolomoitsev@math.uni-goettingen.de}
\thanks{Supported by the German Research Foundation (DFG),
project KO 5804/3-1.}
}

\date{}

\begin{document}

\maketitle

\begin{abstract}
We study sampling recovery in $L_p(\mathbb T^d)$
for intersections of Lipschitz classes defined by smoothness
conditions in two different metrics $L_p$ and $L_q$, where
$1\le p<q\le\infty$.
Our focus is on the low-smoothness range, where the
$L_p$-smoothness does not exceed the classical continuity
threshold $d/p$.
We determine the exact orders of the sampling numbers
for both linear and general recovery algorithms,
including critical smoothness cases in which logarithmic
factors appear or uniform convergence fails.
To analyze sampling approximation in the low-smoothness range,
we establish new inequalities for averaged moduli of
smoothness with different integrability parameters and
derive sharp estimates for these moduli and the corresponding
semi-discrete measures of smoothness.
These estimates provide error bounds for linear sampling
operators acting on functions whose smoothness may lie
below the classical continuity threshold.
Finally, we determine the Kolmogorov widths of the
intersection classes and quantify the gap between
optimal sampling recovery and approximation by
finite-dimensional subspaces. 
\end{abstract}


\section{Introduction}

Sampling recovery is one of the classical problems of approximation
theory and numerical analysis. Given a function $f$ and a finite
set of sampling nodes $X_N=\{\xi_1,\ldots,\xi_N\}\subset\mathbb T^d$,
one seeks to reconstruct $f$ from the data
$f(\xi_1),\ldots,f(\xi_N)$ and to estimate the error of the
corresponding recovery procedure in a prescribed norm.
Two natural target metrics are the uniform norm and the
$L_p$-norm. In the uniform setting, the use of point values
is natural for continuous functions. In $L_p$, however,
point evaluation is not a continuous functional in general.
This raises the question of which smoothness assumptions
allow one to obtain meaningful and sharp $L_p$-error
estimates for recovery from function values.

For classical smoothness spaces, sampling recovery is well understood
when the smoothness is sufficiently high. For example, the embedding
$B^s_{p,\theta}(\mathbb T^d)\hookrightarrow C(\mathbb T^d)$
holds for $s>d/p$, so that point evaluations are well defined
and continuous.
In this range, sharp estimates for sampling numbers
are well known for Besov and Triebel--Lizorkin spaces.
In particular, for suitable parameters,
\begin{equation}
\label{eq:classical-sampling}
g_N(B^s_{p,\theta},L_r)
\asymp
g_N^{\rm lin}(B^s_{p,\theta},L_r)
\asymp
N^{-s/d+(1/p-1/r)_+},
\end{equation}
see \cite{NT06,V07}.
Here $B^s_{p,\theta}$ denotes the unit ball
of the corresponding Besov space and $g_N$ and $g_N^{\rm lin}$ denote the sampling
and linear sampling numbers, respectively.
For the same classes, the Kolmogorov widths in the
$L_p$-metric satisfy
\begin{equation}
\label{eq:classical-widths}
d_N(B^s_{p,\theta},L_p)
\asymp N^{-s/d},
\end{equation}
see, e.g., \cite{V08}.
Thus, for recovery in the $L_p$-metric, sampling numbers
and Kolmogorov widths have the same asymptotic order
in this smoothness range.
However, for Sobolev spaces with smoothness
close to the continuity threshold, logarithmic
gaps between sampling numbers and approximation
numbers may occur even in the Hilbert space
setting; see \cite{HKNV22}.
More general relations between sampling numbers and
Kolmogorov widths have been investigated in
\cite{HNV08,KU20,DKU23};
see also the recent survey \cite{KU26}.

The situation changes considerably when the smoothness
falls below the classical continuity threshold.
In this range, point evaluations are no longer necessarily
controlled by the smoothness norm, and sampling recovery
may fail to converge uniformly.
For certain Sobolev embeddings, this was
established in \cite{H08}.

The limiting case $s=d/p$ requires particular attention.
While $B^{d/p}_{p,1}(\mathbb T^d)$ embeds into
$C(\mathbb T^d)$, the corresponding embedding fails for
$B^{d/p}_{p,\infty}(\mathbb T^d)$.
Moreover, in the one-dimensional case $1<p<\infty$,
Kolyada \cite{K79} showed that even the condition
$\omega_1(f,h)_p=o(h^{1/p})$ does not guarantee
continuity, by constructing bounded functions
that are essentially discontinuous at every point.

Nevertheless, sampling operators can still be defined for functions
with small smoothness and even for discontinuous functions. Indeed,
an operator of the form
\[
G_Nf=\sum_{\xi\in X_N}f(\xi)\varphi_{\xi,N}
\]
is well defined whenever the values $f(\xi)$ are finite and the
reconstruction functions $\varphi_{\xi,N}$ belong to the target space. The main difficulty is therefore not the definition of the operator, but obtaining error estimates in terms of the smoothness of $f$.

One of the classical tools for this purpose is the averaged modulus of smoothness $\tau_r(f,\delta)_p$. Averaged moduli were systematically investigated by Sendov and Popov \cite{SP88} and have been used in one-sided approximation, interpolation, numerical integration, and sampling approximation; see, for example,
\cite{H89,HI90,BBSTV06,BMVSV14,PX94}. In particular, for a wide class of sampling operators and every bounded function $f$, one has a direct estimate of the form
\[
\|f-G_Nf\|_p
\lesssim \tau_r(f,N^{-1/d})_p.
\]

The averaged modulus is related to the ordinary modulus of smoothness
in $L_p(\mathbb T^d)$ by inequalities of the form
\begin{equation}
\label{eq:classical-tau}
\tau_r(f,\delta)_p
\lesssim
\delta^{d/p}
\int_0^\delta
\frac{\omega_r(f,t)_p}{t^{d/p}}
\,\frac{dt}{t},
\end{equation}
see, e.g.,~\cite{H89}. Consequently, if $\omega_r(f,t)_p\lesssim t^\alpha$ with
$\alpha>d/p$, then
$\tau_r(f,\delta)_p\lesssim\delta^\alpha$.
Thus, this classical relation leads to the same threshold $d/p$
that appears in the embedding of Besov spaces into continuous
functions. If $\alpha\le d/p$, however, the above estimate no longer
yields the desired bound for the averaged modulus.

Another difficulty is that direct estimates
through averaged moduli do not generally admit
inverse inequalities of the same type for
sampling operators. This problem was studied
in our previous works \cite{KL23,K26,KBS26},
where sampling-adapted semi-discrete measures
of smoothness were introduced. For example,
in the one-dimensional periodic case, such
a measure can be written in the form
\[
\Omega(f,1/n;X_n)_p
=
\omega_r(f,1/n)_p
+
\|f-f_{1/n}\|_{\ell_p(X_n)},
\]
where $f_{1/n}$ is a Steklov mean of $f$.
The first term describes the ordinary
$L_p$-smoothness of $f$, while the second
measures the difference between the function
and its local averages at the sampling nodes.

Under appropriate assumptions on the sampling
operators, these measures lead to complementary
direct and inverse estimates. In particular, for the
corresponding operators $G_n$ and any finite-valued
function $f\in L_p(\mathbb T)$, one obtains, for every
$\gamma\in(0,2)$, characterizations of approximation
rates of the form
\begin{equation}\label{eq:individual-sampling-equivalence}
\|f-G_nf\|_p=O(n^{-\gamma})
\quad\Longleftrightarrow\quad
\Omega(f,1/n;X_n)_p=O(n^{-\gamma}).
\end{equation}
In contrast to estimates based only on the ordinary
$L_p$-smoothness, this characterization is not restricted
by the threshold $\gamma>1/p$. Thus, the semi-discrete
measures describe the approximation properties of
individual functions also in the low-smoothness range.
On the other hand, the semi-discrete measures
depend explicitly on the sampling nodes.
This raises the question of whether, in the
low-smoothness range, the additional information
about point values can be replaced by smoothness
conditions that are independent of the sampling set.

In this paper, we investigate this question by imposing
smoothness conditions in two different metrics. Alongside
the usual $L_p$-smoothness condition, we impose an additional
condition in $L_q$, where $q>p$.
More precisely, for $1\le p<q\le\infty$,
$\alpha>0$, $\beta\ge0$, and integers
$r>\alpha$ and $s>\beta$, we consider
the intersections of Lipschitz classes
\[
\mathcal F_{p,q}^{\alpha,\beta}(\mathbb T^d)
=
\left\{
 f\in C(\mathbb T^d):
 \|f\|_\infty\le1,
 \ \omega_r(f,t)_p\le t^\alpha,
 \ \omega_s(f,t)_q\le t^\beta
\right\},
\]
where the smoothness conditions are imposed
for all $0<t\le1$.

Our main interest is in the range
$d/q<\beta<\alpha<d/p$, together with
the critical cases $\alpha=d/p$ and
$\beta=d/q$. These classes allow us to
investigate how additional smoothness
in the stronger $L_q$-metric affects
sampling recovery when the $L_p$-smoothness
does not exceed the continuity threshold.

In the interior range
$d/q<\beta<\alpha<d/p$,
it is instructive to compare our results
with what can be obtained by interpolation. The two smoothness conditions imply embeddings into
intermediate Besov spaces $B^\sigma_{v,\infty}$, where
\[
\frac1v=\frac{1-\theta}{p}+\frac{\theta}{q},
\qquad
\sigma=(1-\theta)\alpha+\theta\beta.
\]
The critical relation $\sigma=d/v$ is attained for
\[
\rho=
\frac{\beta d/p-\alpha d/q}
{\beta-\alpha+d/p-d/q}.
\]
Thus, interpolation places
$\mathcal F_{p,q}^{\alpha,\beta}$
in $B^\rho_{v,\infty}$ with $\rho=d/v$,
exactly at the critical continuity threshold.
Choosing a nearby interpolation parameter
with $\sigma>d/v$ and applying
\eqref{eq:classical-sampling} with target metric
$L_p$, we obtain sampling rates of order
$N^{-(\rho-\varepsilon)/d}$ for every
sufficiently small $\varepsilon>0$.
However, this argument does not provide
the limiting rate $N^{-\rho/d}$,
which we show to be attainable and optimal.

Our main result, Theorem~\ref{thm:sampling-numbers-d},
establishes the exact orders of both linear
and nonlinear sampling numbers:
\[
g_N(\mathcal F_{p,q}^{\alpha,\beta},L_p)
\asymp
g_N^{\rm lin}(\mathcal F_{p,q}^{\alpha,\beta},L_p)
\asymp
\begin{cases}
N^{-\rho/d},
& d/q<\beta<\alpha<d/p,\\[2mm]
N^{-1/p}\log(N+1),
& d/q<\beta<\alpha=d/p,\\[2mm]
1,
& \beta=d/q,\quad\beta<\alpha\le d/p.
\end{cases}
\]
Thus, the limiting rate suggested by interpolation
is attained in the interior range. Moreover,
nonlinear recovery does not improve the asymptotic
orders, while on the lower critical line
$\beta=d/q$ there is no uniform convergence
over the class.

The proof of the upper sampling estimates is based
on a refinement of the classical inequality
\eqref{eq:classical-tau}. By combining the
$L_p$- and $L_q$-moduli at different scales,
we establish, for $1\le p<q<\infty$ and
$0<\eta\le\delta/2$, an estimate of the form
\[
\begin{aligned}
\tau_r(f,\delta)_p
\lesssim{}&
\delta^{d/p}
\int_\eta^\delta
\frac{\omega_r(f,t)_p}{t^{d/p}}
\,\frac{dt}{t}+
\delta^{d/q}
\int_0^\eta
\frac{\omega_s(f,t)_q}{t^{d/q}}
\,\frac{dt}{t}.
\end{aligned}
\]
For $q=\infty$, the second term is replaced
by $\omega_s(f,\eta)_\infty$.
Balancing the two terms by a suitable
choice of $\eta$, we obtain the desired
upper estimates for the averaged modulus
on the classes $\mathcal F_{p,q}^{\alpha,\beta}$.

We also obtain sharp estimates for the
semi-discrete smoothness measure $\Omega_r$
on the classes $\mathcal F_{p,q}^{\alpha,\beta}$.
Together with a direct approximation theorem
for linear sampling operators, these estimates
yield the corresponding upper bounds for
sampling errors. The matching lower bounds
are established using localized functions
that vanish at prescribed sampling nodes.

A comparison of the sampling numbers with Kolmogorov widths and recovery from
arbitrary continuous linear information
reveals a substantial difference between
sampling recovery and general approximation.
Indeed, we show in Section~6 that
\[
d_N(\mathcal F_{p,q}^{\alpha,\beta},L_p)
\asymp N^{-\alpha/d},
\]
while recovery from arbitrary continuous
linear information achieves an error of order
at most $N^{-\alpha/d}$.
Consequently, sampling recovery is polynomially
slower in the interior low-smoothness range,
whereas at the upper critical boundary
$\alpha=d/p$ the gap relative to Kolmogorov widths
is only logarithmic.
At the lower critical boundary $\beta=d/q$,
sampling recovery fails to converge uniformly
on the class, although recovery from arbitrary
continuous linear information retains a
polynomial convergence rate.

The paper is organized as follows.
Section~2 contains the necessary notation
and auxiliary results. In Section~3, we
establish mixed-integrability inequalities
for averaged moduli of smoothness and derive
their consequences for the classes under
consideration. Section~4 develops localized
constructions needed for the lower sampling
estimates. Section~5 introduces the
semi-discrete measure of smoothness  and proves
direct estimates for sampling operators,
including their existence on quasi-uniform
sets. Finally, Section~6 establishes the
sharp sampling numbers, derives their
consequences for measures of smoothness and
sampling errors, and compares sampling
recovery with Kolmogorov widths and
arbitrary linear information.

\section{Preliminaries}
\label{sec:preliminaries}

Throughout the paper, $A\lesssim B$ means that $A\le CB$, where
$C>0$ may depend on fixed structural parameters such as $d,p,q,r,s$,
but not on the functions or on the main asymptotic parameters
$\delta,\eta,n,N$. As usual, $A\asymp B$ means that
$A\lesssim B$ and $B\lesssim A$.
For a measurable set $S\subset\mathbb T^d$, let $|S|$ denote
its Lebesgue measure. The convention $1/\infty=0$ is used throughout.

We identify $\mathbb T^d$ with $[0,1)^d$ and use
the maximum norm
\[
|x|=\max_{1\le j\le d}|x_j|.
\]
Distances on $\mathbb T^d$ are understood periodically.
For $x\in\mathbb T^d$ and $t>0$, let
\[
Q(x,t)=\{y\in\mathbb T^d:|x-y|<t\}.
\]

For $r\in\mathbb N$ and $f\in L_p(\T^d)$, $1\le p\le \infty$, the ordinary
modulus of smoothness is defined by
\[
\omega_r(f,t)_p
=
\sup_{|h|\le t}\|\Delta_h^r f\|_{L_p(\mathbb T^d)},
\]
where
\[
\Delta_h^r f(x)
=
\sum_{\nu=0}^r(-1)^{r-\nu}\binom r\nu f(x+\nu h).
\]

The local and averaged moduli of smoothness are defined by
\[
\omega_r(f,x,t)
=
\sup\left\{
|\Delta_h^r f(y)|:
y,y+rh\in x+(-rt,rt)^d
\right\}
\]
and
\[
\tau_r(f,t)_p
=
\|\omega_r(f,\cdot,t)\|_{L_p(\mathbb T^d)}.
\]
These two moduli make sense only for bounded functions $f$.

We shall use the following elementary properties of the ordinary
modulus of smoothness:
\begin{itemize}
\item[$(a)$] $\omega_r(f,t)_p$ is non-negative and non-decreasing in $t$;
\item[$(b)$]
$\omega_r(f+g,t)_p\le\omega_r(f,t)_p+\omega_r(g,t)_p$;
\item[$(c)$]
$\omega_{r+s}(f,t)_p\le2^s\omega_r(f,t)_p$;
\item[$(d)$] for every $\lambda>0$,
\[
\omega_r(f,\lambda t)_p
\le(1+\lambda)^r\omega_r(f,t)_p.
\]
\end{itemize}
For the averaged modulus we use
\begin{itemize}
\item[$(a')$] $\tau_r(f,t)_p$ is non-negative and non-decreasing in $t$;
\item[$(b')$]
$\tau_r(f+g,t)_p\le\tau_r(f,t)_p+\tau_r(g,t)_p$;
\item[$(c')$]
$\tau_{r+1}(f,t)_p\le2\tau_r(f,(r+1)t/r)_p$;
\item[$(d')$] for every fixed $\lambda>0$,
\[
\tau_r(f,\lambda t)_p
\lesssim_\lambda\tau_r(f,t)_p;
\]
\item[$(e')$]
\[
\tau_r(f,t)_p
\lesssim
t^r|f|_{W_p^r}
+
\sum_{\substack{r<|\beta|_1\le d\\ \beta\in\{0,1\}^d}}
t^{|\beta|_1}\|D^\beta f\|_p
\]
for every smooth periodic function $f$.
\end{itemize}
These properties are standard; see, e.g.,
\cite[Ch.~4]{TB04}, \cite[Ch.~2]{DL}, \cite{BMVSV14}, and~\cite{SP88}. For every continuous function $f$ on $\mathbb T^d$,
we also have
\begin{equation}\label{eq:omega-below-tau}
\omega_r(f,t)_p
\le
\tau_r(f,t)_p
\le
\omega_r(f,2t)_\infty.
\end{equation}

We shall consider intersections of Lipschitz classes
defined by two smoothness conditions in different
integrability metrics. For $1\le p<q\le\infty$,
$\alpha>0$, $\beta\ge0$, and fixed integers
$r>\alpha$, $s>\beta$, let
\begin{equation}\label{eq:mixed-class}
\mathcal F_{p,q}^{\alpha,\beta}(\mathbb T^d)
=
\left\{
 f\in C(\mathbb T^d):
 \|f\|_\infty\le1,
 \ \omega_r(f,t)_p\le t^\alpha,
 \ \omega_s(f,t)_q\le t^\beta,
 \ 0<t\le1
\right\}.
\end{equation}
For simplicity, the dependence on $r$ and $s$
is omitted from the notation.

\subsection{Discrete and local estimates for trigonometric polynomials}

Let $\mathcal T_n^d$ denote the trigonometric polynomials
whose spectrum is contained in
$\{k\in\mathbb Z^d:|k|\le n\}$.

We shall use the following upper
Marcinkiewicz--Zygmund inequality for separated
nodes, which follows from
\cite[Theorem~5.5(a)]{FilbirMhaskar2011};
see also \cite{KT26} for a more general framework.

\begin{lemma}\label{lem:irregular-MZ}
Let $1\le p<\infty$ and let $Y_n\subset\Td$
be a finite set satisfying
\[
|x-y|\ge\frac{\gamma}{n},
\qquad x,y\in Y_n,\quad x\ne y,
\]
where $\gamma>0$ is independent of $n$.
Then
\[
n^{-d}\sum_{x\in Y_n}|T(x)|^p
\lesssim \|T\|_p^p,
\qquad T\in\mathcal T_n^d.
\]
\end{lemma}

For $T\in\mathcal T_n^d$, we set
\[
M_tT(x):=\sup_{y\in Q(x,t)}|T(y)|.
\]

\begin{lemma}\label{lem:local-poly}
Let $1\le p\le q\le\infty$ and $T\in\mathcal T_n^d$. Then
\begin{equation}\label{eq:local-maximal}
\|M_tT\|_p
\lesssim
(1+nt)^{d/q}\|T\|_q,
\qquad t>0.
\end{equation}
Consequently,
\begin{equation}\label{eq:tau-poly}
\tau_r(T,t)_p
\lesssim
(1+nt)^{d/q}\|T\|_q.
\end{equation}
\end{lemma}

\begin{proof}
The case $q=\infty$ is immediate.  Assume $q<\infty$ and partition
$\mathbb T^d$ into the cubes
\[
Q_k
=
\prod_{i=1}^d
\left[\frac{k_i}{n},\frac{k_i+1}{n}\right),
\qquad
k\in\{0,\ldots,n-1\}^d.
\]
Put
\[
A_k=\sup_{x\in Q_k}|T(x)|,
\qquad
\Lambda(k)=\{\nu:\operatorname{dist}(Q_\nu,Q_k)\le t\}.
\]
Then
\[
|\Lambda(k)|\lesssim(1+nt)^d,
\qquad
|\{k:\nu\in\Lambda(k)\}|\lesssim(1+nt)^d.
\]  Moreover, for $x\in Q_k$,
\[
M_tT(x)
\le
\max_{\nu\in\Lambda(k)}A_\nu
\le
\left(\sum_{\nu\in\Lambda(k)}A_\nu^q\right)^{1/q}.
\]

Choose $y_k\in\overline{Q_k}$ such that $|T(y_k)|=A_k$. Splitting the
cubes into a fixed finite number of subfamilies, the points $y_k$ in
each subfamily are $c/n$-separated, where $c>0$ depends only on the
dimension. Applying Lemma~\ref{lem:irregular-MZ}
to each subfamily, we obtain
\begin{equation}\label{eq:MZ-maxima}
\frac1{n^d}\sum_kA_k^q\lesssim\|T\|_q^q.
\end{equation}
Since $p/q\le1$,
\[
\begin{aligned}
\|M_tT\|_p^p
&\le
\frac1{n^d}\sum_k
\left(\sum_{\nu\in\Lambda(k)}A_\nu^q\right)^{p/q}\\
&\le
\left(
\frac1{n^d}\sum_k\sum_{\nu\in\Lambda(k)}A_\nu^q
\right)^{p/q}\\
&\lesssim
(1+nt)^{dp/q}
\left(\frac1{n^d}\sum_\nu A_\nu^q\right)^{p/q}\\
&\lesssim
(1+nt)^{dp/q}\|T\|_q^p.
\end{aligned}
\]
This proves \eqref{eq:local-maximal}.

If $y,y+rh\in x+(-rt,rt)^d$, then
$y+\nu h\in x+(-rt,rt)^d$ for
$\nu=0,\ldots,r$. Hence
\[
|\Delta_h^rT(y)|
\le2^r\sup_{z\in Q(x,rt)}|T(z)|,
\]
so
\[
\omega_r(T,x,t)\le2^rM_{rt}T(x).
\]
Taking the $L_p$ norm and applying \eqref{eq:local-maximal} proves
\eqref{eq:tau-poly}.
\end{proof}

\subsection{Dyadic decomposition}

Let $V_j$ be the tensor-product de la Vall\'ee Poussin operators
\[
V_j f(x)
=
\sum_{k\in\mathbb Z^d}v_d(2^{-j}k)\widehat f_k e^{2\pi i k\cdot x},
\qquad
v_d(\xi)=\prod_{\ell=1}^d v_1(\xi_\ell),
\]
where
\[
v_1(t)=
\begin{cases}
1,&|t|\le1,\\
2-|t|,&1<|t|<2,\\
0,&|t|\ge2.
\end{cases}
\]
Set
\[
\Delta_0f=V_0f,
\qquad
\Delta_jf=V_jf-V_{j-1}f,
\quad j\ge1.
\]
We use the standard realization formula, see, e.g.,~\cite{HI90} or~\cite[Ch.~4]{KT21},
\begin{equation}\label{eq:realization}
\|f-V_jf\|_u+2^{-mj}|V_jf|_{W_u^m}
\asymp
\omega_m(f,2^{-j})_u,
\qquad 1\le u\le\infty.
\end{equation}

\begin{lemma}\label{lem:dyadic}
Let $r\in\mathbb N$, $j\ge1$, and $1\le u\le\infty$.
Then
\begin{equation}\label{eq:block-Jackson}
\|\Delta_jf\|_u
\lesssim
\omega_r(f,2^{-j})_u
\end{equation}
and
\begin{equation}\label{eq:low-realization}
\tau_r(V_jf,2^{-j})_u
\lesssim
\omega_r(f,2^{-j})_u.
\end{equation}
\end{lemma}

\begin{proof}
By \eqref{eq:realization} and property~(d),
\[
\begin{aligned}
\|\Delta_jf\|_u
&\le
\|f-V_jf\|_u+\|f-V_{j-1}f\|_u\lesssim
\omega_r(f,2^{-j})_u,
\end{aligned}
\]
which proves \eqref{eq:block-Jackson}.

To establish \eqref{eq:low-realization}, we apply
property~$(e')$ with $t=2^{-j}$.
For $\beta\in\{0,1\}^d$ with $k=|\beta|_1>r$,
choose $\alpha\le\beta$ such that $|\alpha|_1=r$.
By Bernstein's inequality,
\[
\|D^\beta V_jf\|_u
\lesssim
2^{j(k-r)}\|D^\alpha V_jf\|_u.
\]
Consequently, property~$(e')$ and
\eqref{eq:realization} yield
\[
\tau_r(V_jf,2^{-j})_u
\lesssim
2^{-jr}|V_jf|_{W_u^r}
\lesssim
\omega_r(f,2^{-j})_u.
\]
\end{proof}

\section{Mixed-integrability estimates for the \texorpdfstring{$\tau$}{tau}-modulus}
\label{sec:mixed-tau}

The next result refines the classical estimate \(\eqref{eq:classical-tau}\), which uses only the \(L_p\)-smoothness of \(f\) and therefore becomes ineffective at the critical threshold \(d/p\). To overcome this limitation, we incorporate additional \(L_q\)-smoothness, which makes it possible to obtain useful estimates beyond the \(L_p\)-threshold.

\begin{theorem}\label{thm:mixed-tau}
Let $1\le p<q\le\infty$, $r,s\in\mathbb N$, $f\in C(\mathbb T^d)$,
and $0<\eta\le\delta/2$. If $q<\infty$, then
\begin{equation}\label{eq:mixed-tau}
\tau_r(f,\delta)_p
\lesssim
\delta^{d/p}\int_\eta^\delta
\frac{\omega_r(f,t)_p}{t^{d/p+1}}\,dt
+
\delta^{d/q}\int_0^\eta
\frac{\omega_s(f,t)_q}{t^{d/q+1}}\,dt.
\end{equation}
If $q=\infty$, then
\begin{equation}\label{eq:mixed-tau-infty}
\tau_r(f,\delta)_p
\lesssim
\delta^{d/p}\int_\eta^\delta
\frac{\omega_r(f,t)_p}{t^{d/p+1}}\,dt
+
\omega_s(f,\eta)_\infty.
\end{equation}
\end{theorem}

\begin{proof}
We first consider the case $q<\infty$.
Choose integers $J<M$ such that
\[
2^{-J}\asymp\delta,
\qquad
2^{-M}\asymp\eta,
\]
and decompose $f$ into three frequency components:
\[
f=V_Jf+\sum_{j=J+1}^{M}\Delta_jf
+\sum_{j>M}\Delta_jf.
\]
By the subadditivity property~$(b')$ of the averaged modulus,
\[
\tau_r(f,\delta)_p
\le
\tau_r(V_Jf,\delta)_p
+
\sum_{j=J+1}^{M}\tau_r(\Delta_jf,\delta)_p
+
\sum_{j>M}\tau_r(\Delta_jf,\delta)_p.
\]

For the low-frequency component, Lemma~\ref{lem:dyadic}
and the scaling properties~$(d')$ and $(d)$ yield
\[
\tau_r(V_Jf,\delta)_p
\lesssim
\omega_r(f,\delta)_p.
\]
For $J<j\le M$, by Lemmas~\ref{lem:local-poly}
and~\ref{lem:dyadic},  we have
\[
\tau_r(\Delta_jf,\delta)_p
\lesssim
(\delta2^j)^{d/p}\omega_r(f,2^{-j})_p.
\]
Summing over $J<j\le M$ and passing from the dyadic
sum to an integral, we obtain
\begin{equation}\label{eq:mixed-tau-middle}
\sum_{j=J+1}^{M}\tau_r(\Delta_jf,\delta)_p
\lesssim
\delta^{d/p}\int_\eta^\delta
\frac{\omega_r(f,t)_p}{t^{d/p+1}}\,dt.
\end{equation}

For the high-frequency tail, we use the
additional $L_q$-smoothness of $f$ to obtain,
for $j>M$,
\[
\tau_r(\Delta_jf,\delta)_p
\lesssim
(\delta2^j)^{d/q}\omega_s(f,2^{-j})_q,
\]
which, after summation, gives
\begin{equation}\label{eq:mixed-tau-tail}
\sum_{j>M}\tau_r(\Delta_jf,\delta)_p
\lesssim
\delta^{d/q}\int_0^\eta
\frac{\omega_s(f,t)_q}{t^{d/q+1}}\,dt.
\end{equation}

It remains to incorporate the low-frequency contribution
into the first integral. Since $\eta\le\delta/2$,
properties~$(a)$ and~$(d)$ imply
\[
\begin{aligned}
\delta^{d/p}\int_\eta^\delta
\frac{\omega_r(f,t)_p}{t^{d/p+1}}\,dt
&\ge
\delta^{d/p}\int_{\delta/2}^\delta
\frac{\omega_r(f,t)_p}{t^{d/p+1}}\,dt\\
&\gtrsim
\omega_r(f,\delta/2)_p
\gtrsim
\omega_r(f,\delta)_p.
\end{aligned}
\]
Combining this observation with
\eqref{eq:mixed-tau-middle} and
\eqref{eq:mixed-tau-tail} proves
\eqref{eq:mixed-tau}.

For $q=\infty$, we use the decomposition
\[
f=V_Jf+\sum_{j=J+1}^{M}\Delta_jf+(f-V_Mf).
\]
The low- and intermediate-frequency components
are treated as above. The remaining term can be
estimated directly in the uniform norm:
\[
\begin{aligned}
\tau_r(f-V_Mf,\delta)_p
&\le 2^r\|f-V_Mf\|_\infty\\
&\lesssim \omega_s(f,2^{-M})_\infty\lesssim \omega_s(f,\eta)_\infty.
\end{aligned}
\]
Together with \eqref{eq:mixed-tau-middle} and
the estimate for the low-frequency component,
this yields \eqref{eq:mixed-tau-infty}.
\end{proof}

\begin{corollary}\label{cor:power}
Let $1\le p<q\le\infty$, $f\in C(\mathbb T^d)$,
and suppose that
\[
\omega_r(f,t)_p=O(t^\alpha),
\qquad
\omega_s(f,t)_q=O(t^\beta).
\]
\begin{enumerate}%[label={\rm(\roman*)}]
\renewcommand{\labelenumi}{\rm(\roman{enumi})}
\item If $d/q<\beta<\alpha<d/p$, then
\[
\tau_r(f,\delta)_p=O(\delta^\rho).
\]

\item If $d/q<\beta<\alpha=d/p$, then
\[
\tau_r(f,\delta)_p
=
O\left(\delta^{d/p}\log\frac1\delta\right).
\]
\end{enumerate}
\end{corollary}

\begin{proof}
In case (i), Theorem~\ref{thm:mixed-tau} gives
\[
\tau_r(f,\delta)_p
=
O\left(
\delta^{d/p}\eta^{\alpha-d/p}
+
\delta^{d/q}\eta^{\beta-d/q}
\right).
\]
Set
\[
\lambda=
\frac{d/p-d/q}
{\beta-\alpha+d/p-d/q}>1.
\]
Then $\rho$ and $\lambda$ satisfy the balance identities
\begin{equation}\label{eq:balance-identities}
\rho-\frac dp+\lambda\left(\frac dp-\alpha\right)=0,
\qquad
\rho-\frac dq+\lambda\left(\frac dq-\beta\right)=0.
\end{equation}
Taking $\eta\asymp\delta^\lambda$ and using
\eqref{eq:balance-identities} gives
$\tau_r(f,\delta)_p=O(\delta^\rho)$.

In case (ii), Theorem~\ref{thm:mixed-tau} yields
\[
\tau_r(f,\delta)_p
=
O\left(
\delta^{d/p}\log\frac{\delta}{\eta}
+
\delta^{d/q}\eta^{\beta-d/q}
\right).
\]
Taking
$\eta
\asymp
\delta^{(d/p-d/q)/(\beta-d/q)}$
gives the required estimate.
\end{proof}

We now show that the averaged modulus still has a nontrivial power decay when the additional \(L_q\)-smoothness approaches the lower critical threshold. In particular, for \(q=\infty\), even a very weak quantitative continuity condition is sufficient to obtain a power rate. Some such additional control is nevertheless essential. As will be shown later in Proposition~\ref{prop:Omega-arbitrarily-slow}, without it the averaged modulus may tend to zero arbitrarily slowly.

\begin{proposition}
\label{prop:log-refinements}
Let $f\in C(\mathbb T^d)$, $1\le p<q\le\infty$, and
$\omega_r(f,t)_p=O(t^{d/p}).$
\begin{enumerate}%[label={\rm(\roman*)}]
\renewcommand{\labelenumi}{\rm(\roman{enumi})}
\item If $q<\infty$, $\gamma>1$, and
\[
\omega_s(f,t)_q
=
O\left(t^{d/q}(\log(e/t))^{-\gamma}\right),
\]
then
\[
\tau_r(f,\delta)_p
=
O\left(
\delta^{((\gamma-1)d/p+d/q)/\gamma}
\right).
\]

\item If $q=\infty$, $\gamma>0$, and
\[
\omega_s(f,t)_\infty
=
O\left((\log(e/t))^{-\gamma}\right),
\]
then
\[
\tau_r(f,\delta)_p
=
O\left(
\delta^{\gamma d/[p(\gamma+1)]}
\right).
\]
\end{enumerate}
\end{proposition}

\begin{proof}
For (i), Theorem~\ref{thm:mixed-tau} yields
\[
\tau_r(f,\delta)_p
=
O\left(
\delta^{d/p}\log\frac e\eta
+
\delta^{d/q}(\log(e/\eta))^{1-\gamma}
\right).
\]
Taking
\[
\log\frac e\eta
\asymp
\delta^{-(d/p-d/q)/\gamma}
\]
gives the required estimate.

For (ii),
\[
\tau_r(f,\delta)_p
=
O\left(
\delta^{d/p}\log\frac e\eta
+
(\log(e/\eta))^{-\gamma}
\right).
\]
Taking
\[
\log\frac e\eta
\asymp
\delta^{-d/[p(\gamma+1)]}
\]
gives the required estimate.
\end{proof}

The case $p=1$ has an additional one-dimensional feature. Recall that, when $d=1$, the first-order critical condition $\omega_1(f,t)_1=O(t)$ is equivalent to $f$ being equal almost everywhere to a function of bounded variation; see, e.g.,~\cite[p.~53]{DL}. This condition is stronger than the higher-order critical assumption $\omega_r(f,t)_1=O(t)$, $r>1$, used above and in the definition of the corresponding Lipschitz classes. This special endpoint case leads to a refinement of Theorem~\ref{thm:mixed-tau}, in which the first-order modulus plays a distinguished role.

\begin{proposition}
\label{prop:p1-endpoint}
Let $d=1$, $1<q\le\infty$, $r,s\in\mathbb N$,
$f\in C(\mathbb T)$, and $0<\eta\le\delta/2$. Then
\[
\tau_r(f,\delta)_1
\lesssim
\frac{\delta}{\eta}\omega_1(f,\eta)_1
+
\begin{cases}
\displaystyle
\delta^{1/q}\int_0^\eta
\frac{\omega_s(f,t)_q}{t^{1/q+1}}\,dt,
&q<\infty,\\[3mm]
\omega_s(f,\eta)_\infty,
&q=\infty.
\end{cases}
\]
\end{proposition}

\begin{proof}
Choose $M$ so that
$2^{-M-1}<\eta\le2^{-M}$. By $(b')$, we have
\begin{equation}\label{eq:p1-decomposition}
\tau_r(f,\delta)_1
\le \tau_r(V_Mf,\delta)_1+\tau_r(f-V_Mf,\delta)_1.
\end{equation}
For the first term, properties~$(c')$ and~$(e')$,
together with the realization formula~\eqref{eq:realization},
yield
\begin{equation*}
\tau_r(V_Mf,\delta)_1
\lesssim
\tau_1(V_Mf,r\delta)_1
\lesssim
\delta\,\|(V_Mf)'\|_1
\lesssim
\frac{\delta}{\eta}\omega_1(f,\eta)_1.
\end{equation*}
Assume now that $q<\infty$. Using property~$(b')$, Lemmas
\ref{lem:local-poly} and \ref{lem:dyadic}, and the fact that
$\delta2^j\gtrsim1$ for $j>M$, we obtain
\[
\begin{aligned}
\tau_r(f-V_Mf,\delta)_1
&\le
\sum_{j>M}\tau_r(\Delta_jf,\delta)_1\\
&\lesssim
\delta^{1/q}\sum_{j>M}2^{j/q}\omega_s(f,2^{-j})_q\\
&\lesssim
\delta^{1/q}\int_0^\eta
\frac{\omega_s(f,t)_q}{t^{1/q+1}}\,dt.
\end{aligned}
\]
If $q=\infty$, then the realization formula~\eqref{eq:realization} gives instead
\[
\tau_r(f-V_Mf,\delta)_1
\le
2^r\|f-V_Mf\|_\infty
\lesssim
\omega_s(f,\eta)_\infty.
\]
Combining these estimates with
\eqref{eq:p1-decomposition} completes the proof.
\end{proof}

\section{Localized functions and smoothness estimates}
\label{sec:sharpness}

To establish the lower bounds for sampling recovery
in Section~6, we construct functions that vanish
at arbitrary prescribed sampling nodes while retaining
a sufficiently large $L_p$-norm. The following lemma
provides the basic localized functions used in these
constructions.

\begin{lemma}\label{lem:localized-functions}
Let $X\subset\mathbb T^d$ be finite, $0<\eta\le1/4$, and
$m\in\mathbb N$.  There exists $H_{\eta,X}\in C^\infty(\mathbb T^d)$
such that
\[
0\le H_{\eta,X}\le1,
\qquad
H_{\eta,X}(x)=1\quad\text{for}\quad x\in X,
\]
and
\begin{equation}\label{eq:H-support}
\supp H_{\eta,X}
\subset
\bigcup_{x\in X} Q(x,3\eta/8),
\qquad
|\supp H_{\eta,X}|
\le
\left(\frac34\right)^d|X|\eta^d.
\end{equation}
Moreover, for $1\le u\le\infty$,
\begin{equation}\label{eq:H-difference}
\|\Delta_h^mH_{\eta,X}\|_u
\lesssim
(|X|\eta^d)^{1/u}
\min\left\{1,\left(\frac{|h|}{\eta}\right)^m\right\}.
\end{equation}
\end{lemma}

\begin{proof}
Choose a non-negative function $\psi\in C_c^\infty(\mathbb R^d)$ such
that
\[
\int_{\mathbb R^d}\psi=1,
\qquad
\supp\psi\subset Q(0,1/8),
\]
and put
\[
U_{\eta,X}=\bigcup_{x\in X}Q(x,\eta/4),
\qquad
\psi_\eta(x)=\eta^{-d}\psi(x/\eta),
\qquad
H_{\eta,X}=\mathbf1_{U_{\eta,X}}*\psi_\eta,
\]
where the convolution is taken on $\mathbb T^d$ and $\psi_\eta$
denotes the periodic extension of $\eta^{-d}\psi(\cdot/\eta)$.
It is clear that $0\le H_{\eta,X}\le1$. If $x\in X$ and
$y\in\supp\psi_\eta$, then
\[
x-y\in Q(x,\eta/4)\subset U_{\eta,X}.
\]
Therefore
\[
H_{\eta,X}(x)
=
\int_{\mathbb T^d}\mathbf1_{U_{\eta,X}}(x-y)\psi_\eta(y)\,dy=\int_{\mathbb T^d}\psi_\eta(y)\,dy
=1.
\]
Moreover, by the definition of $H_{\eta,X}$ and the support property of $\psi_\eta,$
\[
\supp H_{\eta,X}
\subset
\bigcup_{x\in X}Q(x,3\eta/8),
\]
which implies \eqref{eq:H-support}.

To prove \eqref{eq:H-difference}, we first estimate the
derivatives of $H_{\eta,X}$. For every multi-index
$\nu$ with $|\nu|_1=m$, we have
\[
D^\nu H_{\eta,X}
=
\mathbf1_{U_{\eta,X}}*D^\nu\psi_\eta,
\qquad
\|D^\nu H_{\eta,X}\|_\infty\lesssim\eta^{-m}.
\]
Since these derivatives are supported in
$\supp H_{\eta,X}$, estimate~\eqref{eq:H-support}
implies
\[
\|D^\nu H_{\eta,X}\|_u
\lesssim
(|X|\eta^d)^{1/u}\eta^{-m}.
\]
Using the integral representation of finite
differences for $|h|<\eta$ and the elementary
estimate for $|h|\ge\eta$, we obtain
\[
\begin{aligned}
\|\Delta_h^m H_{\eta,X}\|_u
&\lesssim
\begin{cases}
|h|^m\displaystyle\sum_{|\nu|_1=m}
\|D^\nu H_{\eta,X}\|_u,
& |h|<\eta,\\[2mm]
\|H_{\eta,X}\|_u,
& |h|\ge\eta
\end{cases}\\[2mm]
&\lesssim
(|X|\eta^d)^{1/u}
\min\left\{1,\left(\frac{|h|}{\eta}\right)^m\right\},
\end{aligned}
\]
which proves \eqref{eq:H-difference}.
\end{proof}

\begin{lemma}\label{lem:multilevel-concentration}
Let $X\subset\Td$ be finite, $a>0$,
$0<2^M\eta\le1/4$, $m\in\N$, and $0<\sigma<m$.
For the function
\[
F=a\sum_{\nu=0}^M H_{2^\nu\eta,X},
\]
the following estimates hold.
If $1\le u<\infty$, then
\begin{equation}\label{eq:multilevel-modulus}
\sup_{0<t\le1}t^{-\sigma}\omega_m(F,t)_u
\lesssim
a|X|^{1/u}
\begin{cases}
(2^M\eta)^{d/u-\sigma},
& \sigma<d/u,\\[1mm]
1,
& \sigma=d/u,\\[1mm]
\eta^{d/u-\sigma},
& \sigma>d/u.
\end{cases}
\end{equation}
For $u=\infty$, we have
\begin{equation}\label{eq:multilevel-modulus-infty}
\sup_{0<t\le1}t^{-\sigma}
\omega_m(F,t)_\infty
\lesssim a\eta^{-\sigma},
\end{equation}
while
\[
\sup_{0<t\le1}\omega_m(F,t)_\infty
\lesssim a(M+1).
\]
\end{lemma}

\begin{proof}
Applying Lemma~\ref{lem:localized-functions}
to each summand of $F$, we obtain
\begin{equation}\label{eq:multilevel-start}
\omega_m(F,t)_u
\lesssim
a|X|^{1/u}
\sum_{\nu=0}^M
(2^\nu\eta)^{d/u}
\min\left\{1,
\left(\frac{t}{2^\nu\eta}\right)^m
\right\}.
\end{equation}
Since $m>\sigma$, we have
\[
\min\{1,z^m\}\le z^\sigma,
\qquad z>0,
\]
and therefore
\[
t^{-\sigma}\omega_m(F,t)_u
\lesssim
a|X|^{1/u}
\sum_{\nu=0}^M
(2^\nu\eta)^{d/u-\sigma}.
\]
Evaluating the geometric sum gives
\eqref{eq:multilevel-modulus} for
$1\le u<\infty$ and $\sigma\ne d/u$,
as well as \eqref{eq:multilevel-modulus-infty}
for $u=\infty$.

It remains to consider the critical case
$\sigma=d/u$ for $1\le u<\infty$.
Splitting the sum in \eqref{eq:multilevel-start}
according to $2^\nu\eta\le t$ and
$2^\nu\eta>t$, and using $m>d/u$, we obtain
\[
\begin{aligned}
\omega_m(F,t)_u
&\lesssim
a|X|^{1/u}
\left(
\sum_{2^\nu\eta\le t}
(2^\nu\eta)^{d/u}
+
t^m
\sum_{2^\nu\eta>t}
(2^\nu\eta)^{d/u-m}
\right)\\
&\lesssim
a|X|^{1/u}t^{d/u}.
\end{aligned}
\]
Finally,
\[
\omega_m(F,t)_\infty
\le 2^m\|F\|_\infty
\le 2^m a(M+1),
\]
which proves the last assertion.
\end{proof}

\begin{lemma}\label{lem:membership}
Let $1\le p<q\le\infty$,
$0<\alpha\le d/p$, $\beta\ge d/q$,
$r,s\in\N$ with $r>\alpha$ and $s>\beta$,
and let $X\subset\Td$ be finite.
For $a>0$ and
$0<2^M\eta\le1/4$, consider
\[
F=a\sum_{\nu=0}^M H_{2^\nu\eta,X}.
\]
Suppose that
\begin{equation}\label{eq:membership-conditions}
a|X|^{1/p}(2^M\eta)^{d/p-\alpha}\lesssim1,
\qquad
a|X|^{1/q}\eta^{d/q-\beta}\lesssim1.
\end{equation}
If $q=\infty$ and $\beta=0$, assume in addition
that $a(M+1)\lesssim1$.
Then
\begin{equation}\label{eq:membership-conclusion}
\omega_r(F,t)_p\lesssim t^\alpha,
\qquad
\omega_s(F,t)_q\lesssim t^\beta,
\qquad 0<t\le1.
\end{equation}
The same estimates hold for $F+c$,
$c\in\mathbb R$.
\end{lemma}

\begin{proof}
Applying Lemma~\ref{lem:multilevel-concentration}
with $(u,m,\sigma)=(p,r,\alpha)$ and using
the first condition in \eqref{eq:membership-conditions},
we obtain
\[
\omega_r(F,t)_p\lesssim t^\alpha.
\]
Similarly, for $q<\infty$, applying the same lemma
with $(u,m,\sigma)=(q,s,\beta)$ and using
the second condition yields
\[
\omega_s(F,t)_q\lesssim t^\beta.
\]
For $q=\infty$, the latter estimate follows
from \eqref{eq:multilevel-modulus-infty} when
$\beta>0$, and from the last assertion of
Lemma~\ref{lem:multilevel-concentration}
when $\beta=0$, using the additional assumption
$a(M+1)\lesssim1$.

Since finite differences annihilate constants,
the estimates also hold for $F+c$.
\end{proof}

\section{Semi-discrete smoothness and sampling operators}
\label{sec:semi-discrete-sampling}

Let $X_n\subset\mathbb T^d$ be finite sets satisfying
\begin{equation}\label{eq:quasi-uniform-simple}
|X_n|\asymp n^d,
\qquad
|x-y|\ge\frac{\gamma}{n},
\quad x,y\in X_n,\quad x\ne y,
\end{equation}
where $\gamma>0$ is fixed. Set
\[
\|g\|_{\ell_p(X_n)}
=
\left(n^{-d}\sum_{x\in X_n}|g(x)|^p\right)^{1/p}.
\]

\subsection{Semi-discrete modulus of smoothness}

Fix \(r\in\mathbb N\) and put \(k=\lceil r/2\rceil\). For \(0<h<1/2\), define
\[
A_hf(x)=\frac{1}{(2h)^d}\int_{[-h,h]^d}f(x+u)\,du,
\]
and
\[
f_{h,k}(x)=
-\frac{2}{\binom{2k}{k}}
\sum_{j=1}^k(-1)^j
\binom{2k}{k-j}A_{jh/k}f(x).
\]
This is a generalized Steklov average with an integrable kernel supported in \(Q(0,h)\). Moreover,
\begin{equation}\label{eq:steklov-difference}
f(x)-f_{h,k}(x)
=
\frac{(-1)^k}{\binom{2k}{k}(2h/k)^d}
\int_{[-h/k,h/k]^d}
\Delta_u^{2k}f(x-ku)\,du.
\end{equation}

We define the semi-discrete measure of smoothness by
\[
\Omega_r(f,1/n;X_n)_p
=
\|f-f_{1/n,k}\|_{\ell_p(X_n)}
+\omega_r(f,1/n)_p.
\]

As discussed in the Introduction, the semi-discrete
measure $\Omega_r$ is better adapted to the approximation
of individual functions by sampling operators than the
averaged modulus; see also \cite{KL23,K26,KBS26}.
We now relate $\Omega_r$ to the averaged modulus
$\tau_r$ in order to transfer the upper estimates
obtained in Section~3 to the semi-discrete setting.
The following lemma provides the required comparison.

\begin{lemma}\label{lem:averaging}
Let $1\le p<\infty$, $r\in\N$, and
$k=\lceil r/2\rceil$.
Suppose that $X_n$ satisfies
\eqref{eq:quasi-uniform-simple}.
Then, for $h=1/n$,
\begin{align}
\|f-f_{h,k}\|_p
&\lesssim \omega_r(f,h)_p,
\label{eq:average-error}\\
\|f_{h,k}\|_{\ell_p(X_n)}
&\lesssim \|f\|_p.
\label{eq:average-discrete}
\end{align}
Both inequalities hold for arbitrary
$f\in L_p(\Td)$.
Moreover, for $f\in C(\Td)$,
\begin{equation}\label{eq:Omega-tau-comparison}
\Omega_r(f,1/n;X_n)_p
\lesssim \tau_r(f,1/n)_p.
\end{equation}
\end{lemma}

\begin{proof}
Estimate~\eqref{eq:average-error} follows from
\eqref{eq:steklov-difference}, Minkowski's inequality,
and the standard properties of the ordinary modulus
of smoothness. The kernel representation of $f_{h,k}$ and the
bounded overlap of the cubes $Q(x,h)$,
$x\in X_n$, yield
\[
\begin{aligned}
\|f_{h,k}\|_{\ell_p(X_n)}^p
&\lesssim
\sum_{x\in X_n}\int_{Q(x,h)}|f(y)|^p\,dy\lesssim \|f\|_p^p,
\end{aligned}
\]
which proves \eqref{eq:average-discrete}.

Finally, \eqref{eq:Omega-tau-comparison} follows
by adapting the argument of
\cite[Proposition~2.1]{K26} to periodic cubes.
Indeed, \eqref{eq:steklov-difference},
the bounded overlap of the cubes $Q(x,Ch)$,
$x\in X_n$, and the standard properties
of averaged moduli give
\[
\|f-f_{h,k}\|_{\ell_p(X_n)}
\lesssim \tau_r(f,h)_p.
\]
Together with
$\omega_r(f,h)_p\lesssim\tau_r(f,h)_p$,
this proves \eqref{eq:Omega-tau-comparison}.
\end{proof}

\subsection{Linear sampling operators}

Consider the linear sampling operators
\begin{equation}\label{eq:sampling-operator-d}
G_nf
=
\sum_{x\in X_n}f(x)\varphi_{x,n}.
\end{equation}

We assume the following stability and approximation
properties, which hold uniformly with respect to $n$
and the sampling nodes:
\begin{align}
\|G_nf\|_p
&\lesssim
\|f\|_{\ell_p(X_n)},
\tag{$A_1$}\label{eq:sampling-A1-d}\\
\|T-G_nT\|_p
&\lesssim
n^{-r}|T|_{W_p^r},
\qquad T\in\mathcal T_n^d.
\tag{$A_2$}\label{eq:sampling-A3-d}
\end{align}
These properties are the multidimensional periodic
counterparts of those used in \cite{KL23}.

\begin{theorem}\label{thm:sampling-direct-d}
Let $1\le p<\infty$ and $r\in\N$.
Suppose that $X_n$ satisfies
\eqref{eq:quasi-uniform-simple} and that
$G_n$ satisfies \eqref{eq:sampling-A1-d}
and \eqref{eq:sampling-A3-d}.
Then, for every $f\in L_p(\Td)$
with prescribed point values at the sampling nodes,
\begin{equation}\label{eq:sampling-direct-Omega-d}
\|f-G_nf\|_p
\lesssim
\Omega_r(f,1/n;X_n)_p.
\end{equation}
In particular, for $f\in C(\Td)$,
\begin{equation}\label{eq:sampling-direct-tau-d}
\|f-G_nf\|_p
\lesssim
\Omega_r(f,1/n;X_n)_p
\lesssim
\tau_r(f,1/n)_p.
\end{equation}
\end{theorem}

\begin{proof}
The argument is analogous to those of
\cite[Theorem~3.1]{K26} and the corresponding
periodic result in \cite{KL23}.
Put $h=1/n$ and choose $j$ such that
$2^{j+1}\le n$ and $2^j\asymp n$.
Consider the trigonometric polynomial
$P=V_jf\in\mathcal T_n^d$.
By the realization formula \eqref{eq:realization},
\[
\|f-P\|_p+n^{-r}|P|_{W_p^r}
\lesssim \omega_r(f,h)_p.
\]
Furthermore, Lemmas~\ref{lem:averaging}
and~\ref{lem:irregular-MZ} yield
\[
\begin{aligned}
\|f-P\|_{\ell_p(X_n)}
&\le \|f-f_{h,k}\|_{\ell_p(X_n)}
 +\|(f-P)_{h,k}\|_{\ell_p(X_n)}
 +\|P_{h,k}-P\|_{\ell_p(X_n)}\\
&\lesssim
\|f-f_{h,k}\|_{\ell_p(X_n)}
+\|f-P\|_p+\|P_{h,k}-P\|_p\\
&\lesssim
\|f-f_{h,k}\|_{\ell_p(X_n)}
+\omega_r(f,h)_p.
\end{aligned}
\]
Consequently, assumptions
\eqref{eq:sampling-A1-d} and
\eqref{eq:sampling-A3-d} imply
\[
\begin{aligned}
\|f-G_nf\|_p
&\le \|f-P\|_p+\|P-G_nP\|_p
+\|G_n(f-P)\|_p\\
&\lesssim
\|f-P\|_p
+n^{-r}|P|_{W_p^r}
+\|f-P\|_{\ell_p(X_n)}\\
&\lesssim
\omega_r(f,h)_p
+\|f-f_{h,k}\|_{\ell_p(X_n)}\\
&=\Omega_r(f,h;X_n)_p.
\end{aligned}
\]

The last assertion follows from
\eqref{eq:Omega-tau-comparison}.
\end{proof}

\begin{corollary}\label{cor:sampling-rates-d}
Let $1\le p<q\le\infty$, and suppose that
$X_n$ satisfies \eqref{eq:quasi-uniform-simple}.
Let $G_n$ satisfy \eqref{eq:sampling-A1-d}
and \eqref{eq:sampling-A3-d} with $r>\alpha$.
Then
\[
\sup_{f\in\mathcal F_{p,q}^{\alpha,\beta}}
\|f-G_nf\|_p
\lesssim
\begin{cases}
n^{-\rho},
& d/q<\beta<\alpha<d/p,\\[2mm]
n^{-d/p}\log(n+1),
& d/q<\beta<\alpha=d/p,\\[2mm]
1,
& \beta=d/q,\quad \beta<\alpha\le d/p.
\end{cases}
\]
\end{corollary}

\begin{proof}
The first two estimates follow from
Theorem~\ref{thm:sampling-direct-d}
and Corollary~\ref{cor:power}.
The last estimate follows from
$\tau_r(f,1/n)_p\lesssim\|f\|_\infty\le1$.
\end{proof}

\subsection{Existence of sampling operators on quasi-uniform sets}

We now impose the additional covering condition
\begin{equation}\label{eq:sampling-covering-d}
\sup_{z\in\Td}\inf_{x\in X_n}|z-x|
\le \frac{H}{n},
\end{equation}
where $H>0$ is independent of $n$.

\begin{proposition}\label{prop:quasi-uniform-operators}
Let $1\le p<\infty$, and suppose that $X_n$ satisfies
\eqref{eq:quasi-uniform-simple} and
\eqref{eq:sampling-covering-d}.
Then, for every fixed $r\in\N$, there exist
linear sampling operators $G_n$ of the form
\eqref{eq:sampling-operator-d} satisfying
\eqref{eq:sampling-A1-d} and
\eqref{eq:sampling-A3-d}
for all sufficiently large $n$, with constants
independent of $n$.
\end{proposition}

\begin{proof}
By the quadrature theorem in
\cite[Theorem~5.7]{FilbirMhaskar2011},
the conditions \eqref{eq:quasi-uniform-simple}
and \eqref{eq:sampling-covering-d} imply that,
for some sufficiently small $\theta>0$,
there exist nonnegative weights
$w_{x,n}\lesssim n^{-d}$ such that
\begin{equation}\label{eq:quasi-uniform-cubature}
\int_{\Td}P(z)\,dz
=
\sum_{x\in X_n}w_{x,n}P(x),
\qquad P\in\mathcal T_{\lfloor\theta n\rfloor}^d.
\end{equation}
Here the bound on the weights follows by applying
the quadrature formula to a nonnegative
tensor-product Fej\'er kernel centered at each node.

Choose $j$ such that $2^j\asymp n$,
with a sufficiently small implied constant
to ensure $3\cdot2^j\le\theta n$.
Let $K_j$ be the convolution kernel of
the de la Vall\'ee Poussin operator $V_j$
introduced in Section~2. Then
\begin{equation}\label{eq:Vj-properties}
K_j\in\mathcal T_{2^{j+1}}^d,
\qquad
\|K_j\|_1\lesssim1,
\qquad
V_jP=P,\quad P\in\mathcal T_{2^j}^d.
\end{equation}
Define
\[
G_nf(z)=
\sum_{x\in X_n}w_{x,n}f(x)K_j(z-x).
\]
Since \eqref{eq:quasi-uniform-cubature}
is exact on $\mathcal T_{3\cdot2^j}^d$,
we have
\begin{equation}\label{eq:Gn-reproduction}
G_nP=P,\qquad P\in\mathcal T_{2^j}^d.
\end{equation}

To prove \eqref{eq:sampling-A1-d},
let $1\le p<\infty$ and let $p'$ be its
conjugate exponent. By \eqref{eq:Vj-properties},
duality, H\"older's inequality,
Lemma~\ref{lem:irregular-MZ},
and Young's convolution inequality,
\[
\begin{aligned}
\|G_nf\|_p
&=
\sup_{\|g\|_{p'}\le1}
\left|\int_{\Td}G_nf(z)\overline{g(z)}\,dz\right|\\
&\le
\sup_{\|g\|_{p'}\le1}
\sum_{x\in X_n}
w_{x,n}|f(x)|\,|V_jg(x)|\\
&\lesssim
\|f\|_{\ell_p(X_n)}
\sup_{\|g\|_{p'}\le1}
\|V_jg\|_{\ell_{p'}(X_n)}\\
&\lesssim
\|f\|_{\ell_p(X_n)}
\sup_{\|g\|_{p'}\le1}\|V_jg\|_{p'}\\
&\lesssim \|f\|_{\ell_p(X_n)}.
\end{aligned}
\]
For $p=1$, the corresponding discrete
$\ell_\infty$ estimate is trivial.

Finally, for $T\in\mathcal T_n^d$,
Jackson's inequality gives a polynomial
$P\in\mathcal T_{2^j}^d$ satisfying
\[
\|T-P\|_p
\lesssim 2^{-jr}|T|_{W_p^r}
\lesssim n^{-r}|T|_{W_p^r}.
\]
Therefore, by \eqref{eq:Gn-reproduction},
\eqref{eq:sampling-A1-d}, and
Lemma~\ref{lem:irregular-MZ},
\[
\begin{aligned}
\|T-G_nT\|_p
&\le \|T-P\|_p+\|G_n(T-P)\|_p\\
&\lesssim \|T-P\|_p
\lesssim n^{-r}|T|_{W_p^r},
\end{aligned}
\]
which proves \eqref{eq:sampling-A3-d}.
\end{proof}

\section{Sampling numbers and optimal recovery}
\label{sec:sampling-numbers}

Let $\mathcal F\subset C(\mathbb T^d)$.
For $N\in\mathbb N$, define
\begin{equation}\label{eq:sampling-number-d}
g_N(\mathcal F,L_p)
=
\inf_{x_1,\ldots,x_N\in\mathbb T^d}
\inf_{\Phi:\mathbb C^N\to L_p(\mathbb T^d)}
\sup_{f\in\mathcal F}
\left\|f-\Phi(f(x_1),\ldots,f(x_N))\right\|_p
\end{equation}
and
\begin{equation}\label{eq:linear-sampling-number-d}
g_N^{\rm lin}(\mathcal F,L_p)
=
\inf_{x_1,\ldots,x_N\in\mathbb T^d}
\inf_{\varphi_1,\ldots,\varphi_N\in L_p}
\sup_{f\in\mathcal F}
\left\|f-\sum_{j=1}^Nf(x_j)\varphi_j\right\|_p.
\end{equation}
Clearly,
$g_N(\mathcal F,L_p)\le
g_N^{\rm lin}(\mathcal F,L_p)$.

We now determine the exact asymptotic orders of
the linear and nonlinear sampling numbers for
the classes $\mathcal F_{p,q}^{\alpha,\beta}$.

\begin{theorem}
\label{thm:sampling-numbers-d}
Let $1\le p<q\le\infty$. Then
\begin{enumerate}%[label={\rm(\roman*)}]
\renewcommand{\labelenumi}{\rm(\roman{enumi})}
\item if $d/q<\beta<\alpha<d/p$, then
\begin{equation}\label{eq:sampling-numbers-interior-d}
g_N^{\rm lin}(\mathcal F_{p,q}^{\alpha,\beta},L_p)
\asymp
N^{-\rho/d};
\end{equation}
\item if $\alpha=d/p$ and $d/q<\beta<d/p$, then
\begin{equation}\label{eq:sampling-numbers-critical-d}
g_N^{\rm lin}(\mathcal F_{p,q}^{d/p,\beta},L_p)
\asymp
N^{-1/p}\log(N+1);
\end{equation}
\item if $\beta=d/q$ and $\beta<\alpha\le d/p$, then
\begin{equation}\label{eq:sampling-numbers-lower-critical-d}
g_N^{\rm lin}(\mathcal F_{p,q}^{\alpha,d/q},L_p)
\asymp1.
\end{equation}
\end{enumerate}
The same asymptotic estimates hold for
$g_N(\mathcal F_{p,q}^{\alpha,\beta},L_p)$.
\end{theorem}

\begin{proof}
We first establish the upper bounds.
Choosing uniform grids $X_n$ with $n^d$ nodes,
Proposition~\ref{prop:quasi-uniform-operators}
and Corollary~\ref{cor:sampling-rates-d}
give
\[
g_{n^d}^{\rm lin}
(\mathcal F_{p,q}^{\alpha,\beta},L_p)
\lesssim
\begin{cases}
n^{-\rho},
& d/q<\beta<\alpha<d/p,\\[2mm]
n^{-d/p}\log(n+1),
& d/q<\beta<\alpha=d/p,\\[2mm]
1,
& \beta=d/q,\quad \beta<\alpha\le d/p.
\end{cases}
\]
For arbitrary $N\in\N$, taking
$n=\lfloor N^{1/d}\rfloor$ and using
the monotonicity of sampling numbers
yields the asserted upper bounds
for $g_N^{\rm lin}$ and hence for $g_N$.

For the lower bounds, fix arbitrary sampling nodes
$X=\{x_1,\ldots,x_N\}\subset\Td$
and a reconstruction map
$\Phi:\mathbb C^N\to L_p(\Td)$.
Suppose that $f_X\in\mathcal F_{p,q}^{\alpha,\beta}$
vanishes on $X$.
Since $\mathcal F_{p,q}^{\alpha,\beta}$
is centrally symmetric, the triangle inequality gives
\begin{equation}\label{eq:zero-information-lower}
\max\{\|f_X-\Phi(0)\|_p,
\|-f_X-\Phi(0)\|_p\}
\ge\|f_X\|_p.
\end{equation}
Consequently, any lower bound for $\|f_X\|_p$
uniform over all sampling nodes yields the same
lower bound for $g_N$ and hence for $g_N^{\rm lin}$.

We now construct such functions $f_X$ for the three
parameter ranges considered in the theorem,
using the localized functions $H_{2^\nu\eta,X}$
from Lemma~\ref{lem:localized-functions}.
Throughout the constructions, $c_0>0$ denotes
a sufficiently small constant independent of
$N$ and the sampling nodes.

\noindent\textit{(i)}
Let
\[
\kappa=
\frac{1/p-1/q}{\beta-\alpha+d/p-d/q},
\qquad
\eta=N^{-\kappa},
\qquad
a=c_0N^{-\rho/d}.
\]
Since $\alpha>\beta$, we have $d\kappa>1$.
Moreover, the balance identities
\eqref{eq:balance-identities}, with
$\lambda=d\kappa$, imply
\[
aN^{1/p}\eta^{d/p-\alpha}\asymp c_0,
\qquad
aN^{1/q}\eta^{d/q-\beta}\asymp c_0.
\]

Define
\[
f_X=a(1-H_{\eta,X}).
\]
Then $f_X$ vanishes on $X$ and
$\|f_X\|_\infty\le a\le1$ for sufficiently
large $N$. Since finite differences
annihilate constants, Lemma~\ref{lem:membership}
with $M=0$ shows that
$f_X\in\mathcal F_{p,q}^{\alpha,\beta}$
provided $c_0$ is sufficiently small.

On the other hand, by \eqref{eq:H-support},
\[
|\supp H_{\eta,X}|
\le \left(\frac34\right)^dN\eta^d
\lesssim N^{1-d\kappa}\to0.
\]
Consequently, $f_X=a$ on a set whose
measure tends to one, so that
\[
\|f_X\|_p\gtrsim a
\asymp N^{-\rho/d}.
\]
Together with \eqref{eq:zero-information-lower},
this proves (i).

\noindent\textit{(ii)}
Let
\[
\kappa=
\frac{1/p-1/q}{\beta-d/q}>\frac1d,
\qquad
\eta=N^{-\kappa},
\qquad
a=c_0N^{-1/p}.
\]
For sufficiently large $N$, choose $M\in\mathbb N_0$
such that
\[
\frac{N^{-1/d}}6
<2^M\eta\le
\frac{N^{-1/d}}3.
\]
Since $\kappa>1/d$, such an $M$ exists
and satisfies $M+1\asymp\log N$.

Consider
\[
F_X=a\sum_{\nu=0}^M H_{2^\nu\eta,X},
\qquad
b_N=a(M+1),
\qquad
f_X=b_N-F_X.
\]
Since $H_{2^\nu\eta,X}(x_j)=1$ for every
$x_j\in X$, the function $f_X$ vanishes on $X$.

The choice of parameters gives
\[
aN^{1/p}=c_0,
\qquad
aN^{1/q}\eta^{d/q-\beta}=c_0.
\]
Moreover, $0\le F_X\le b_N$ and
$b_N\asymp c_0N^{-1/p}\log N\to0$.
Thus $\|f_X\|_\infty\le1$ for sufficiently
large $N$. Since finite differences
annihilate constants, Lemma~\ref{lem:membership},
with $c_0$ chosen sufficiently small,
implies that
$f_X\in\mathcal F_{p,q}^{d/p,\beta}$.

By \eqref{eq:H-support},
\[
\supp F_X
\subset
\bigcup_{j=1}^N
Q\left(x_j,\frac38\,2^M\eta\right),
\]
and therefore
\[
|\supp F_X|
\le N\left(\frac34\,2^M\eta\right)^d
\le4^{-d}<\frac12.
\]
Consequently, $f_X=b_N$ on a set
of measure at least $1/2$, which yields
\[
\|f_X\|_p
\gtrsim b_N
\asymp N^{-1/p}\log N.
\]
Together with \eqref{eq:zero-information-lower},
this proves (ii).

\noindent\textit{(iii)}
Choose
\[
M+1=\lceil N^{1/p}\rceil,
\qquad
a=\frac{c_0}{M+1},
\qquad
\eta=\frac{N^{-1/d}}{3\cdot2^M}.
\]
Define
\[
F_X=a\sum_{\nu=0}^M H_{2^\nu\eta,X},
\qquad
f_X=c_0-F_X.
\]
Then $f_X$ vanishes on $X$, and
$0\le f_X\le c_0\le1$.

For every $d/q<\alpha\le d/p$,
the choice of parameters gives
\[
\begin{aligned}
aN^{1/p}(2^M\eta)^{d/p-\alpha}
&\lesssim c_0N^{(\alpha-d/p)/d}
\le c_0,\\
aN^{1/q}\eta^{d/q-\beta}
&=aN^{1/q}\le c_0,
\end{aligned}
\]
where $\beta=d/q$.
If $q=\infty$, we also have
$a(M+1)=c_0$.
Thus, by Lemma~\ref{lem:membership},
$f_X\in\mathcal F_{p,q}^{\alpha,d/q}$
provided $c_0$ is sufficiently small.

Furthermore, \eqref{eq:H-support} yields
\[
\supp F_X
\subset
\bigcup_{j=1}^N Q(x_j,N^{-1/d}/8),
\]
and hence
\[
|\supp F_X|
\le
N\left(\frac{N^{-1/d}}4\right)^d
=4^{-d}<\frac12.
\]
Consequently, $f_X=c_0$ on a set
of measure at least $1/2$, so that
\[
\|f_X\|_p\gtrsim c_0\asymp1.
\]
Together with \eqref{eq:zero-information-lower},
this proves (iii).
\end{proof}

The sharp sampling estimates also determine
the exact worst-case orders of the semi-discrete
measure $\Omega_r$, the averaged modulus $\tau_r$,
and the errors of the corresponding sampling operators.
In particular, the following corollary shows that
$\Omega_r$ and $\tau_r$ have the same worst-case
asymptotic order on the classes
$\mathcal F_{p,q}^{\alpha,\beta}$,
although $\Omega_r$ is more directly related
to sampling recovery; cf.~\eqref{eq:individual-sampling-equivalence}.

\begin{corollary}\label{cor:Omega-tau-sharp-d}
Let $1\le p<q\le\infty$,
$d/q\le\beta<\alpha\le d/p$,
and $r\in\N$ with $r>\alpha$.
Suppose that $X_n$ satisfies
\eqref{eq:quasi-uniform-simple} and
\eqref{eq:sampling-covering-d}.
Then
\[
\sup_{f\in\mathcal F_{p,q}^{\alpha,\beta}}
\Omega_r(f,1/n;X_n)_p
\asymp
\sup_{f\in\mathcal F_{p,q}^{\alpha,\beta}}
\tau_r(f,1/n)_p
\asymp
\begin{cases}
n^{-\rho},
& d/q<\beta<\alpha<d/p,\\[2mm]
n^{-d/p}\log(n+1),
& d/q<\beta<\alpha=d/p,\\[2mm]
1,
& \beta=d/q,\quad\beta<\alpha\le d/p.
\end{cases}
\]
Moreover, the same asymptotic orders hold for
\[
\sup_{f\in\mathcal F_{p,q}^{\alpha,\beta}}
\|f-G_nf\|_p
\]
for every family of operators $G_n$ satisfying
\eqref{eq:sampling-A1-d} and
\eqref{eq:sampling-A3-d}.
\end{corollary}

\begin{proof}
The upper bounds for $\tau_r$ follow from
Corollary~\ref{cor:power} when $\beta>d/q$,
and from the uniform boundedness of
$\mathcal F_{p,q}^{\alpha,d/q}$ when $\beta=d/q$.

By Proposition~\ref{prop:quasi-uniform-operators},
the conditions on $X_n$ ensure the existence
of sampling operators satisfying
\eqref{eq:sampling-A1-d} and
\eqref{eq:sampling-A3-d}.
For any such operators $G_n$,
Theorem~\ref{thm:sampling-direct-d}
and \eqref{eq:Omega-tau-comparison} give
\[
\begin{aligned}
g_{|X_n|}(\mathcal F_{p,q}^{\alpha,\beta},L_p)
&\le
\sup_{f\in\mathcal F_{p,q}^{\alpha,\beta}}
\|f-G_nf\|_p\\
&\lesssim
\sup_{f\in\mathcal F_{p,q}^{\alpha,\beta}}
\Omega_r(f,1/n;X_n)_p\\
&\lesssim
\sup_{f\in\mathcal F_{p,q}^{\alpha,\beta}}
\tau_r(f,1/n)_p.
\end{aligned}
\]
Since $|X_n|\asymp n^d$, Theorem~
\ref{thm:sampling-numbers-d} provides
the corresponding lower bounds for
the first term, and hence for all
three quantities.
Together with the upper bounds for $\tau_r$,
this proves the asserted asymptotic estimates.
\end{proof}

On the lower critical line $\beta=d/q$,
the worst-case errors remain bounded away
from zero. Nevertheless, for each fixed
continuous function, the corresponding
sampling errors and measures of smoothness
tend to zero. The following result shows
that, despite this convergence, no universal
rate can be prescribed.

\begin{proposition}\label{prop:Omega-arbitrarily-slow}
Let $1\le p<q\le\infty$, $\beta=d/q$, and
$d/q<\alpha\le d/p$.
Suppose that $X_n$ satisfies
\eqref{eq:quasi-uniform-simple} and
\eqref{eq:sampling-covering-d}, and that
$G_n$ satisfies \eqref{eq:sampling-A1-d}
and \eqref{eq:sampling-A3-d}, with $r>\alpha$.
For every function $\Psi:(0,1)\to(0,\infty)$
such that $\Psi(t)\to0$ as $t\to0+$,
there exists $f\in\mathcal F_{p,q}^{\alpha,d/q}$
such that
\[
\limsup_{n\to\infty}
\frac{\|f-G_nf\|_p}{\Psi(1/n)}
=\infty.
\]
The same assertion holds with the numerator
replaced by $\Omega_r(f,1/n;X_n)_p$
or $\tau_r(f,1/n)_p$.
\end{proposition}

\begin{proof}
Consider the Banach space $E$ of continuous
functions equipped with the norm
\[
\|f\|_E
=
\max\left\{
\|f\|_\infty,\,
\sup_{0<t\le1}t^{-\alpha}\omega_r(f,t)_p,\,
\sup_{0<t\le1}t^{-d/q}\omega_s(f,t)_q
\right\}.
\]
Then $\mathcal F_{p,q}^{\alpha,d/q}$ is
the closed unit ball of $E$.

By Theorem~\ref{thm:sampling-numbers-d}(iii)
and $|X_n|\asymp n^d$,
\[
\begin{aligned}
\|I-G_n\|_{E\to L_p}
&=
\sup_{f\in\mathcal F_{p,q}^{\alpha,d/q}}
\|f-G_nf\|_p\\
&\ge
g_{|X_n|}
(\mathcal F_{p,q}^{\alpha,d/q},L_p)
\gtrsim1.
\end{aligned}
\]
Consequently,
\[
\sup_{n\ge2}
\left\|
\frac{I-G_n}{\Psi(1/n)}
\right\|_{E\to L_p}
=\infty.
\]
By the Banach--Steinhaus theorem, there exists
$f\in E$ such that
\[
\limsup_{n\to\infty}
\frac{\|f-G_nf\|_p}{\Psi(1/n)}
=\infty.
\]
After multiplying $f$ by a positive constant,
we may assume that
$f\in\mathcal F_{p,q}^{\alpha,d/q}$.

Finally, Theorem~\ref{thm:sampling-direct-d}
and \eqref{eq:Omega-tau-comparison} give
\[
\|f-G_nf\|_p
\lesssim
\Omega_r(f,1/n;X_n)_p
\lesssim
\tau_r(f,1/n)_p,
\]
which proves the remaining assertions.
\end{proof}

\begin{remark}
Although the convergence on the lower critical line
$\beta=d/q$ can be arbitrarily slow, even a mild
logarithmic improvement in smoothness can lead
to a polynomial rate of convergence.
Indeed, Proposition~\ref{prop:log-refinements},
together with \eqref{eq:Omega-tau-comparison}
and Theorem~\ref{thm:sampling-direct-d},
provides such rates for $\tau_r$, $\Omega_r$,
and the errors of sampling operators satisfying
\eqref{eq:sampling-A1-d} and
\eqref{eq:sampling-A3-d}.
\end{remark}

\subsection{Comparison with Kolmogorov widths and linear information}

We conclude this section by comparing the sampling
numbers with Kolmogorov widths and recovery from
arbitrary continuous linear information.

For a bounded set $\mathcal F\subset L_p(\Td)$,
recall that
\[
d_N(\mathcal F,L_p)
=
\inf_{\substack{V\subset L_p(\Td)\\ \dim V<N}}
\sup_{f\in\mathcal F}
\inf_{g\in V}\|f-g\|_p.
\]

The asymptotic behavior of Kolmogorov widths
and sampling numbers for classical smoothness
classes has been extensively investigated;
see, e.g., \cite{NT06,V07,V08,HKNV22}.
The following proposition quantifies the loss
caused by restricting recovery to point evaluations
for the low-regularity classes considered here.

\begin{proposition}\label{prop:kolmogorov-widths}
Let $1\le p<q\le\infty$ and
$d/q\le\beta<\alpha\le d/p$. Then
\[
d_N(\mathcal F_{p,q}^{\alpha,\beta},L_p)
\asymp N^{-\alpha/d}.
\]
Consequently,
\[
\frac{g_N(\mathcal F_{p,q}^{\alpha,\beta},L_p)}
{d_N(\mathcal F_{p,q}^{\alpha,\beta},L_p)}
\asymp
\begin{cases}
N^{(\alpha-\rho)/d},
& d/q<\beta<\alpha<d/p,\\[2mm]
\log(N+1),
& d/q<\beta<\alpha=d/p,\\[2mm]
N^{\alpha/d},
& \beta=d/q,\quad\beta<\alpha\le d/p.
\end{cases}
\]
\end{proposition}

\begin{proof}
Let $U^\alpha_{v,\infty}$ denote the unit ball
of $B^\alpha_{v,\infty}(\Td)$.
The standard Besov embeddings give
\[
c\,U^\alpha_{\infty,\infty}
\subset
\mathcal F_{p,q}^{\alpha,\beta}
\subset
C\,U^\alpha_{p,\infty}
\]
for suitable constants $c,C>0$.
Hence, by the monotonicity of Kolmogorov widths
and the classical estimates for periodic Besov
classes (see, e.g., \cite[Theorem~E]{RomanyukEtAl2023}),
\[
d_N(\mathcal F_{p,q}^{\alpha,\beta},L_p)
\asymp N^{-\alpha/d}.
\]
The remaining assertions follow from
Theorem~\ref{thm:sampling-numbers-d}.
\end{proof}

\begin{remark}
1) A similar gap arises when sampling recovery
is compared with recovery from arbitrary
continuous linear information.
Let $e_N^{\rm all}(\mathcal F,L_p)$ denote
the optimal worst-case recovery error based
on $N$ arbitrary continuous linear functionals
on $L_p(\Td)$, allowing nonlinear reconstruction.
Standard trigonometric approximation using
$O(N)$ Fourier coefficients gives
\[
e_N^{\rm all}
(\mathcal F_{p,q}^{\alpha,\beta},L_p)
\lesssim N^{-\alpha/d}.
\]
Thus, restricting linear information to point
evaluations leads to at least the same asymptotic
loss as that exhibited by the Kolmogorov widths.

2) An interesting feature is the upper critical
case $\alpha=d/p$: for every
$d/q<\beta<d/p$, restricting recovery
to point evaluations causes only a
logarithmic loss, independently of the
dimension $d$.
\end{remark}

\end{document}